\documentclass[a4paper,12pt]{article}
\usepackage[english]{babel}
\usepackage[T2A]{fontenc}
\usepackage[cp1251]{inputenc}
\usepackage[english]{babel}
\usepackage{amsthm}
\usepackage[tbtags]{amsmath}
\usepackage{amsfonts,amssymb}
\begin{document}

\newtheorem{theorem}{Theorem}
\newtheorem{lemma}{Lemma}
\newtheorem{proposition}{Proposition}
\newtheorem{Cor}{Corollary}

\begin{center}
{\large\bf Centrally Essential Torsion-Free Rings of Finite Rank}
\end{center}
\begin{center}
O.V. Lyubimtsev\footnote{Nizhny Novgorod State University; email: oleg lyubimcev@mail.ru .}, A.A. Tuganbaev\footnote{National Research University `MPEI', Lomonosov Moscow State University; email: tuganbaev@gmail.com .}
\end{center}

\textbf{Key words:} centrally essential ring, quasi-invariant ring, faithful Abelian group.

\textbf{Abstract.}
It is proved that centrally essential rings, whose additive groups of finite rank are torsion-free groups of finite rank, are quasi-invariant but not necessarily invariant.
Torsion-free Abelian groups of finite rank with centrally essential endomorphism rings are faithful.

The work of O.V. Lyubimtsev is done under the state assignment No~0729-2020-0055. A.A. Tuganbaev is supported by Russian Scientific Foundation, project 16-11-10013P.

\textbf{MSC2010 database 16R99; 20K30}

\section{Introduction}

We consider only associative rings with non-zero identity elements. For a ring $R$, we denote by $C(R)$, $J(R)$ and $N(R)$ the center, the Jacobson radical and the upper radical of the ring $R$, respectively. For Abelian groups, we use the additive notation.

\textbf{1.1. Centrally essential rings.}\\
A ring $R$ is said to be \textsf{centrally essential} if 
for every non-zero element $a\in R$, there exist non-zero central elements $x,y$ with $ax=y$. \footnote{It is clear that a ring $R$ with center $C$ is centrally essential if and only if the module $R_{C}$ is an essential extension of the module $C_{C}$.} 

Centrally essential rings with non-zero identity elements are studied in \cite{MT18}, \cite{MT19}, \cite{MT19b}, \cite{MT19c}, \cite{MT20}, \cite{MT20b}, \cite{MT20c}. Every centrally essential semiprime ring with $1\ne 0$ is commutative; see \cite[Proposition 3.3]{MT18}. 
Examples of non-commutative group algebras over fields are given in \cite{MT18}. 
In addition, the Grassman algebra over a three-dimensional vector space over the field of order 3 also is a finite non-commutative centrally essential ring; see \cite{MT19}.
In \cite{MT19c}, there is an example of a centrally essential ring $R$ whose factor ring with respect to the prime radical of $R$ is not a PI ring.
Abelian groups with centrally essential endomorphism rings are considered in \cite{LT20}.

\textbf{1.2. Torsion-free rings of finite rank and faithful Abelian groups.}\\
A ring $R$ is called a \textsf{torsion-free ring of finite rank} (a  \textsf{tffr ring}) if the additive group $(R, +)$  of $R$ is an Abelian torsion-free group of finite rank.\\ A ring $R$ is said to be \textsf{right invariant} (resp., \textsf{right quasi-invariant}) if every right ideal (resp., maximal right ideal) of $R$ is an ideal of $R$.
The words such as "an invariant ring" (resp., "a quasi-invariant ring") mean "a right and left invariant ring" (resp., "right and left quasi-invariant ring").\\
An Abelian group $A$ is said to be \textsf{faithful} if $IA\neq A$ for every proper right ideal $I$ of the endomorphism ring $\text{End}\,A$ of the group $A$.

It is known that torsion-free groups of finite rank with commutative endomorphism rings are faithful (see \cite[Theorem 5.9]{Ar82}). Faticoni got a stronger result, he proved that torsion-free groups of finite rank with right invariant endomorphism rings also are faithful groups (see \cite[Lemma 3.1]{Fat88}). Detailed information on faithful Abelian groups is contained in \cite[Chapter VI, \S 33, 34]{KMT03}.

The main results of this paper are Theorems 1.3 and 1.4.

\textbf{1.3. Theorem.} Every centrally essential torsion-free ring of finite rank is a quasi-invariant ring which is not necessarily right or left invariant\footnote{In Example 2.4 of Section 2 of this paper, we give an example of a centrally essential torsion-free ring of finite rank which is not a right or left invariant ring.}.

\textbf{1.4. Theorem.} An Abelian torsion-free group of finite rank with centrally essential endomorphism ring is a faithful group.

We give some definitions used in the paper. A right ideal $I$ of the ring $R$ is said to be \textsf{essential} if $I\cap J\neq 0$ for every non-zero right ideal $J$ in $R$. In this case, one says that $R$ is an \textsf{essential extension} of $I$.

A right ideal $B$ is said to be \textsf{closed} if it coincides with any right ideal which is an essential extension of $B$. 

For a right ideal $I$ of the ring $R$, every right ideal $J$ of $R$,  which is maximal with respect to the property $I\cap J = 0$, is called a \textsf{$\cap$-complement} for $I$.
It is easy to verify that $I\bigoplus J$ is an essential right ideal in $R$ in this case (e.g., see \cite{G76}).

An element $r$ of the ring $R$ is called \textsf{right regular} or a \textsf{left non-zero-divisor} if the relation $rx = 0$ implies  the relation $x = 0$ for every $x\in R$.

A subgroup $B$ of the Abelian group $A$ is said to be \textsf{pure} if the equation $nx = b\in B$, which is solvable in the group $A$, is solvable in $B$. 
 
\section{Properties of Ideals of Centrally Essential Rings} 

\textbf{2.1. Proposition.} Let $R$ be a centrally essential ring with center $C(R)$. If there exists a maximal right ideal $M$ of the ring $R$ which is not an ideal, then $C(R)\cap \left(\bigcap\limits_{n\ge 1}M^n\right)\ne 0$.

\textbf{Proof.} We assume the contrary. Then there exist two elements $m\in M$ and $a\in R$ with $am\notin M$. Since $M$ is a maximal right ideal, there exist elements $b\in R$ and $m'\in M$ such that $1 = amb + m'$. Since $am\notin M$, we have 
$a\ne 0$. Since $R$ is a centrally essential ring, there exist 
non-zero elements $c,d\in C(R)$ such that $ac = d\ne 0$. Then
$$
c = (amb + m')c = (ac)mb + m'c = mbd + m'c\in M,
$$
$(ac)mb\in M^2$  and $m'c\in M^2$. Therefore, $c = (ac)mb + m'c\in M^2$ and $(ac)mb, m'c\in M^3$. Then $c\in M^3$. By repeating a similar argument, we obtain that $0\ne c\in C(R)\cap \left(\bigcap\limits_{n=1}^{+\infty}M^n\right)$.~\hfill$\square$

\textbf{2.2. Proposition.} Let $R$ be a centrally essential ring. Then every minimal right ideal in $R$ is contained in the center $C(R)$ of the ring $R$. In particular, every minimal right ideal of the ring $R$ is an ideal in $R$.

\textbf{Proof.} Let $I$ be a non-zero minimal right ideal. Since the ring $R$ is centrally essential, we have that for every non-zero element $a\in I$, there exist two non-zero central elements $c,d\in R$ such that $ac = d\neq 0$. Then $d\in C(R)\cap I$. Consequently, $C(R)\cap I\neq 0$. In addition, $C(R)\cap I = K$ is an ideal in $R$, $K\subseteq I$. Since $I$ is a minimal right ideal, $I = K$.~\hfill$\square$

\textbf{2.3. Proposition.} If a closed right ideal of a centrally essential ring contains a right regular element, then the right ideal is an ideal.

\textbf{Proof.} Let $J$ be a closed right ideal of the centrally essential ring $R$ and $r\in R$. If $0\neq X\le rJ + J$, then
$X\cap C(R)\neq 0$. We take a non-zero element $x\in X\cap C(R)$. We have $x = ra + b$ for some $a, b\in J$. By assumption, $J$ contains a right regular element $y$.
In this case, we have
$$
0\neq yx = yra + yb\in J\cap X.
$$ 
Consequently, the right ideal $rJ + J$ is an essential extension of $J$. Since $J$ is a closed ideal,  $J = rJ + J$ and $J$ is an ideal.~\hfill$\square$

\textbf{2.4. Example.} We consider the subring $\mathcal{R}$ in the ring $M_{7}(R)$ of all $7\times 7$ matrices over the commutative domain $R$ consisting of the matrices $A$ of the form 
$$
A = 
\left(\begin{matrix}
\alpha & a & b & c & d & e & f\\
0 & \alpha & 0 & b & 0 & 0 & d\\
0 & 0 & \alpha & 0 & 0 & 0 & e\\
0 & 0 & 0 & \alpha & 0 & 0 & 0\\
0 & 0 & 0 & 0 & \alpha & 0 & a\\
0 & 0 & 0 & 0 & 0 & \alpha & b\\
0 & 0 & 0 & 0 & 0 & 0 & \alpha\\
\end{matrix}\right).
$$
Let $A'\in \mathcal{R}$ with $a' = a + 1$ and let the remaining components of $A'$ coincide with corresponding components of the matrix $A$. 
Then $AA'\neq A'A$ if $a\neq 0$ and $b\neq 0$. 
Thus, the ring $\mathcal{R}$ is not commutative. It is easy to see that $C(\mathcal{R})$ consist of matrices
$$
C = 
\left(\begin{matrix}
\alpha & 0 & 0 & c & d & e & f\\
0 & \alpha & 0 & 0 & 0 & 0 & d\\
0 & 0 & \alpha & 0 & 0 & 0 & e\\
0 & 0 & 0 & \alpha & 0 & 0 & 0\\
0 & 0 & 0 & 0 & \alpha & 0 & 0 \\
0 & 0 & 0 & 0 & 0 & \alpha & 0 \\
0 & 0 & 0 & 0 & 0 & 0 & \alpha\\
\end{matrix}\right).
$$
Let $A \in \mathcal{R}$ and let $a\neq 0$ or $b\neq 0$. 
We take the matrix $B\in C(\mathcal{R})$ with $d = a$, $e = b$ and zeroes on the remaining positions. Then $0\neq AB\in C(\mathcal{R})$.
Thus, $\mathcal{R}$ is a centrally essential ring.

We consider the right ideal $I$ of $\mathcal{R}$ consisting of the matrices of the form 
$$
B = 
\left(\begin{matrix}
0 & 0 & b & 0 & 0 & 0 & f\\
0 & 0 & 0 & b & 0 & 0 & 0\\
0 & 0 & 0 & 0 & 0 & 0 & 0\\
0 & 0 & 0 & 0 & 0 & 0 & 0\\
0 & 0 & 0 & 0 & 0 & 0 & 0 \\
0 & 0 & 0 & 0 & 0 & 0 & b \\
0 & 0 & 0 & 0 & 0 & 0 & 0\\
\end{matrix}\right).
$$
It is directly verified that $I$ is not an ideal in $\mathcal{R}$. In addition, $I$ is a closed right ideal. Indeed, the ideal of $\mathcal{R}$, which has only $c$ as a non-zero component, is a $\cap$-complement for $I$. 

At the same time, the closed left ideal $J$ of $\mathcal{R}$ consisting of elements are the matrices
$$
D = 
\left(\begin{matrix}
0 & a & 0 & 0 & 0 & 0 & f\\
0 & 0 & 0 & 0 & 0 & 0 & 0\\
0 & 0 & 0 & 0 & 0 & 0 & 0\\
0 & 0 & 0 & 0 & 0 & 0 & 0\\
0 & 0 & 0 & 0 & 0 & 0 & a \\
0 & 0 & 0 & 0 & 0 & 0 & 0 \\
0 & 0 & 0 & 0 & 0 & 0 & 0\\
\end{matrix}\right).
$$
is not an ideal. The ideal which has only $c$ as a non-zero component, is also a $\cap$-complement for $J$. 

The above example can be generalized to the case of matrices of any order $n$. Namely, the subring $\mathcal{R}$ of the ring $M_n(R)$ of all $n\times n$ matrices over the commutative ring $R$ consisting of matrices $A$ of the form
$$
A = 
\left(\begin{matrix}
\alpha & a_{12} & a_{13} & a_{14} & a_{15} & \ldots & a_{1n-2} & a_{1n-1} & a_{1n}\\
0 & \alpha & 0 & a_{13} & 0 & \ldots & 0 & 0 & a_{1n-2}\\
0 & 0 & \alpha & 0 & 0 & \ldots & 0 & 0 & a_{1n-1}\\
0 & 0 & 0 & \alpha & 0 & \ldots & 0 & 0 & 0\\
\ldots & \ldots & \ldots & \ldots & \ldots & \ldots & \ldots & \ldots &\ldots\\
0 & 0 & 0 & 0 & 0 & \ldots & 0 & 0 & 0\\
0 & 0 & 0 & 0 & 0 & \ldots & \alpha & 0 & a_{12}\\
0 & 0 & 0 & 0 & 0 & \ldots & 0 & \alpha & a_{13}\\
0 & 0 & 0 & 0 & 0 & \ldots & 0 & 0 & \alpha\\
\end{matrix}\right),
$$
is non-commutative centrally essential ring which is not invariant. 

\section{The proof of Theorem 1.3} 

For any two elements $a, b$ of the ring $R$, their commutator $ab - ba$ is denoted by $[a, b]$. 
 
\textbf{3.1. Proposition.} Let $R$ be a centrally essential tffr ring. Then the ring $R/N(R)$ is commutative.

\textbf{Proof.} Since $R$ is a centrally essential ring, $\mathbb{Q}R = \mathbb{Q}\otimes R$ also is a centrally essential ring; see \cite[Proposition 3.1]{LT20}. In addition, the ring $\mathbb{Q}R$  is Artinian, since it is a finite-dimensional $\mathbb{Q}$-algebra. It is well known that $N(\mathbb{Q}R) = J(\mathbb{Q}R)$; e.g., see \cite[Proposition 9.1(c)]{Ar82}. Then the ring $\mathbb{Q}R/N(\mathbb{Q}R)$ is commutative \cite[Theorem 1.5]{MT19b}. Let $a, b\in R$. By considering \cite[Proposition 9.1(c)]{Ar82}, we have
$$
[a, b]\in N(\mathbb{Q}R)\cap R = N(R).
$$
Consequently, the ring $R/N(R)$ is commutative.~\hfill$\square$
 
\textbf{3.2. Corollary.} If $R$ is a centrally essential tffr ring, then the ring $R/J(R)$ is commutative.

\textbf{Proof.} Indeed, we have for $a, b\in R$:
$$
[a, b]\in N(R)\subseteq J(R).\eqno \square$$
 
\textbf{3.3. The completion of the proof of Theorem 1.3.}

Let $R$ be a centrally essential ring, $r\in R$, and let $M$ be a maximal right ideal of the ring $R$. We have to prove that $rM\subseteq M$. We assume the contrary. Then $rm\notin M$ for some element $m\in M$. Since $N(R)\subseteq M$, we have $r, m\notin N(R)$. Consequently, $\overline{r} = r + N(R)$ and $\overline{m} = m + N(R)$ are non-zero elements of the commutative ring
$R/N(R)$. Then $\overline{r}\,\overline{m} = \overline{m}\,\overline{r}$; this implies that $[r, m]\in N(R)$. Since $mr\in M$, we have $rm\in M$. This is a contradiction. ~\hfill$\square$

\textbf{3.4. Remark} If we take a tffr ring as the ground ring $R$ in Example 2.4, then the ring $\mathcal{R}$ also is a tffr ring. This implies that a tffr ring is not necessarily right or left invariant.
 
\section{The proof of Theorem 1.4} 

\textbf{4.1. Lemma \cite[Proposition 9.1(c)]{Ar82}} Let $R$ be a tffr ring and $nr\in C(R)$ for some $n\in \mathbb{Z}$, $r\in R$. Then $r\in C(R)$. Thus, $C(R)$ is a pure subgroup of $(R, +)$.
 
\textbf{Proof.} We assume the contrary: $r\notin C(R)$. Then $rx\neq xr$ for some $r\in R$. This implies that $n(rx - xr)\neq 0$ and $(nr)x\neq x(nr)$. This is a contradiction.~\hfill$\square$

\textbf{4.2 Proposition} Let $R$ be a centrally essential tffr ring. Then $R/pR$ is a centrally essential ring.

\textbf{Proof.} Let $r + pR\in R/pR$, where $r\notin pR$. We prove that there exists an element $c'\in C(R)\backslash pR$ such that
$0\neq rc' = d'\in C(R)\backslash pR$. 

Since $R$ is a centrally essential ring, there exists an element $c\in C(R)$ such that $0\neq rc = d\in C(R)$. Let $c\in pR$ and let the $p$-height $h_p(c) = k$,
$c = p^kc'$. It follows from Lemma 4.1 that $c'\in C(R)$. Then $rc = p^krc'\in C(R)$.
We again use Lemma 4.1 and obtain that $0\neq rc' = d'\in C(R)\backslash pR$.~\hfill$\square$

\textbf{4.3. The completion of the proof of Theorem 1.4.} It is well known that for an Abelian torsion-free group $A$ of finite rank, the endomorphism ring $\text{End}\,A$ is a tffr ring. We denote by $R$ the ring $\text{End}\,A$. In \cite{LT20}, it is proved that an Abelian torsion-free group $A$ of finite rank can have non-commutative centrally essential endomorphism ring only if $A$ is a reduced strongly indecomposable group. Following \cite[Lemma 3.1]{Fat88}, we assume that $IA = A$ for some maximal right ideal of the ring $R$. It follows from \cite{Fat88} that $pR\subseteq I$ for some prime integer $p$. By Theorem 1.3, the ring $R$ is quasi-invariant. Therefore, $I$ is a maximal ideal. Since $R/pR$ is a finite ring, the Jacobson radical $J/pR$ of $R$ is nilpotent. We have $J\subseteq I$. Since $R/J$ is a semisimple Artinian ring, $I/J = \overline{e}(R/J)$ for some idempotent $1\neq \overline{e}\in R/J$. Since the ideal $J/pR$ is nilpotent, we can lift $\overline{e}$ to some idempotent $1\neq e\in R/pR$, i.e., $\overline{e} = e + J/pR$. It follows from Proposition 4.2 that the ring $R/pR$ is centrally essential. In a centrally essential ring, all idempotents are central (see \cite[Lemma 2.3]{MT18}). Therefore, the right ideal $e(R/pR)$ is an ideal of $R/pR$. The remaining part of the proof repeat the Faticoni's proof. We give it for completeness.

We remark that 
$$
I/pR = e(R/pR) + J/pR\neq R/pR.
$$
We have
\begin{multline*}
IA/pA = (I/pR)(A/pA) =\\
= [e(R/pR) + J/pR](A/pA) = e(A/pA) + J/pR(A/pA),
\end{multline*}
by our choice of the idempotent $e$. Since $A/pA$ is a finite $R/pR$-module and $I/pR$, $e(R/pR)$ are ideals of the ring $R/pR$, we can apply the Nakayama lemma to the relation $IA/pA = e(A/pA) + J/pR(A/pA)$. As a result, we obtain $IA/pA = e(A/pA)$.
Then the non-zero element $1 - e$ annihilates $IA/pA$. Since $A/pA$ is a faithful left $R/pR$-module, we have $IA/pA\neq A/pA$. Therefore, $IA\neq A$; this is a contradiction. Thus, $A$ is a faithful group.~\hfill$\square$

In Example 2.4, we set $R = \mathbb{Z}$. Then the ring $\mathcal{R}$ is a countable ring whose additive group is a free group of finite rank $n$. Therefore, this ring is the endomorphism ring of some Abelian torsion-free group of rank $n$ (see \cite{Z67}). Therefore, we obtain Corollary 4.4.

\textbf{4.4. Corollary.} For every positive integer $n\ge 7$, there exists a faithful Abelian torsion-free group of rank $n$ with centrally essential endomorphism ring which is not right invariant.

\section{Remarks and Open Questions}

\textbf{5.1. Open question.} Is it true that every essential right ideal of a centrally essential torsion-free ring of finite rank is an ideal?

\textbf{5.2.} Let $A$ be a torsion-free group of finite rank which is a flat module over the endomorphism ring of $A$ (i.e., $A$ is an \textsf{endo-flat} group). In \cite{Fat88}, it is proved that $\text{End}\,A$ is a right invariant ring if and only if $A$ is a faithful group and every $A$-generated subgroup of the group $A$ is fully invariant in $A$. In connection to this fact, we formulate open question 5.3.

\textbf{5.3. Open question.} Let $A$ be an endo-flat torsion-free strongly indecomposable faithful Abelian group of finite rank. Find a group-theoretical property which is is equivalent to the property that the endomorphism ring of $A$ is a centrally essential ring.

\textbf{5.4. Open question.} Is it true that there exists an Abelian group $A$ with centrally essential endomorphism ring $\text{End}\,A$ such that $\text{End}\,A$ is not a ring with polynomial identity?

\end{document}